\documentclass[11pt,twoside,a4paper]{article}
\usepackage{amsmath,amssymb,amsthm}
\usepackage{ifpdf}
\usepackage[shortcuts]{extdash}
\usepackage[utf8]{inputenc}
\usepackage{standalone}
\usepackage{tikz}

\usepackage{microtype}
\usepackage[toc,page]{appendix}

\usepackage{ifthen}
\usepackage{hyperref}

\newtheorem{theorem}{Theorem}

\newtheorem{lemma}{Lemma}

\newtheorem{remark}{Remark}

\newtheorem{proposition}{Proposition}
\newtheorem{corollary}{Corollary}

\DeclareMathOperator{\crrt}{cr}
\DeclareUnicodeCharacter{2212}{ }

\title{Some invariant properties of quasi-Möbius maps}
\author{Loreno Heer}
\date{\today}

\begin{document}
\maketitle

\begin{abstract}
  We investigate properties which remain invariant under the action of quasi-Möbius
  maps of quasi-metric spaces. A metric space is called doubling with constant $D$
  if every ball of finite radius can be covered by at most $D$ balls of half the radius.
  It is shown that the doubling property is an invariant property for
  (quasi\=/)Möbius maps.
  Additionally it is shown that the property of uniform disconnectedness is an
  invariant for (quasi\=/)Möbius maps as well.
\end{abstract}


\medskip


\section{Introduction}

Let $(X, d)$ be a metric space. $X$ is \emph{doubling} if there exists a constant
$D > 0$, such that every ball of finite radius can be covered by at most $D$
balls of half the radius.
$X$ is \emph{uniformly disconnected} if there exists a constant $\theta < 1$,
such that $X$ contains no $\theta$-chain, i.e. a sequence of (at least $3$
distinct) points $(x_0, x_1, \ldots, x_n)$ such that
\[d(x_i, x_{i+1}) \leq \theta d(x_0, x_n).\]
A map $f:(X,d) \to (Y,d')$ is \emph{quasi-Möbius} if it is a homeomorphism and
there exists a homeomorphism $\nu: [0, \infty[ \to [0, \infty[$, such that for
all quadruples $Q = (x_1, x_2, x_3, x_4)$ of distinct points of $X$ and
$Q' := (f(x_1), f(x_2), f(x_3), f(x_4))$, \[\crrt(Q', d') \leq \nu(\crrt(Q, d))\]
holds. Here the \emph{cross-ratio} $\crrt$ is given by
\[\crrt(Q,d) := \frac{d(x_1,x_3)d(x_2,x_4)}{d(x_1,x_4)d(x_2,x_3)}.\]

The aim of this paper is to prove the following two theorems:
\begin{theorem}[Invariance of doubling under quasi-Möbius maps]\label{thm:main}
  Let $(X,d)$ be a doubling space. Let $f :(X,d) \to (Y,d')$ be a quasi-Möbius
  homeomorphism. Then $(Y,d')$ is doubling.
\end{theorem}

\begin{theorem}[Invariance of uniform disconnectedness under quasi-Möbius maps]\label{thm:undisc}
  Let $(X,d)$ be a metric uniformly disconnected space and let $f: (X,d) \to (Y,d')$ be a quasi-Möbius homeomorphism.
  Then $(Y,d')$ is uniformly disconnected.
\end{theorem}

The results are related to results of Lang-Schlichenmaier~\cite{lang_nagata_2005}
and Xie~\cite{xie_nagata_2008} who proved that quasi-symmetric maps respectively
quasi-Möbius maps preserve the Nagata dimension of metric spaces. The present
work has been inspired by the article of Xie~\cite{xie_nagata_2008} and the work
of Väisälä~\cite{vaisala_quasimobius_1984}. We note that a space is doubling if
and only if it has finite Assouad dimension~\cite{mackay_conformal_2010}.
However the Assouad dimension is not a quasi-symmetric (and therefore also not
a quasi-Möbius) invariant~\cite{tyson_lowering_2001}.

We would like to note that we have been informed that \autoref{thm:main}
is also a consequence of a published result of Li-Shanmugalingam \cite{li_preservation_2015}.

It is well known that uniform disconnectedness is invariant under quasi-symmetric
maps~\cite{mackay_conformal_2010, david_fractured_1997}.
However its behaviour under quasi-Möbius maps has not been studied before.

The related property of uniform perfectness has been shown to be invariant under
the metric inversion in \cite{meyer_uniformly_2009}. It is therefore also invariant
under quasi-Möbius maps.

In \autoref{apx:A} we prove a slight generalization of \autoref{thm:main} and \autoref{thm:undisc}
for $K$-quasi-metric spaces.

\section*{Acknowledgments}
The author would like to thank Viktor Schroeder and Urs Lang for several helpful discussions.

\section{Basic Definitions and Preparations}

We introduce the necessary definitions which we will require later.

\subsection{Extended Metrics}

Let $X$ be a set with cardinality at least 3. We call a map $d: X \times X \to [0, \infty]$ an \emph{extended metric} on $X$ if
there exists a set $\Omega(d) \subset X$ with cardinality 0 or 1 and furthermore all of the following requirements are satisfied:
\begin{enumerate}
\item $d_{\vert X \setminus \Omega(d) \times X \setminus \Omega(d)} : X \setminus \Omega(d) \times X \setminus \Omega(d) \to [0, \infty[ $ is a metric;
\item $d(x, \omega) = d(\omega,x) = \infty$ for all $x \in X \setminus \Omega(d)$ and $\omega \in \Omega(d)$;
\item $d(\omega, \omega) = 0$ for $\omega \in \Omega(d)$.
\end{enumerate}

If $\Omega(d)$ is non empty we call $\omega \in \Omega(d)$ the \emph{infinitely remote point} of $X$. By abuse of notation we may
write $\infty$ for the point $\omega$.

\subsection{Doubling Property}
We call a metric space \emph{doubling with constant $D$} if every ball of finite radius can be covered
by at most $D$ balls of half the radius.

\subsection{Uniform Disconnectedness}

For $\theta < 1$ we call a sequence of (at least 3 distinct) points $(x_0, x_1, \ldots, x_n)$ in a metric space $(X,d)$ a \emph{$\theta$-chain} if
\[d(x_i, x_{i+1}) \leq \theta d(x_0,x_n)\]
holds for all $i \in \{0, 1, \ldots, n-1\}$. A metric space is called \emph{uniformly disconnected with constant $\theta$} if
it contains no $\theta$-chains.\footnote{And therefore also no $\theta'$-chains for any $\theta' \leq \theta$.}

\subsection{Quasi-Möbius and Quasi-Symmetric Maps}
%

We call a homeomorphism $f: (X,d) \to (Y,d')$ \emph{$\nu$-quasi-symmetric} if
for all pairwise distinct $x_1,x_2,x_3 \in X$ we have
\[\frac{d'(f(x_1),f(x_2))}{d'(f(x_1),f(x_3))} \leq \nu(\frac{d(x_1,x_2)}{d(x_1,x_3)}).\]
A homeomorphism $f : (X,d) \to (Y, d')$ is called \emph{quasi-symmetric} if it is
$\nu$-quasi-symmetric for some homeomorphism $\nu : [0, \infty[ \to [0, \infty[$.
It is called \emph{symmetric} if for all pairwise distinct $x_1,x_2,x_3 \in X$ we have
\[\frac{d'(f(x_1),f(x_2))}{d'(f(x_1),f(x_3))} = \frac{d(x_1,x_2)}{d(x_1,x_3)}.\]

\section{Invariance of Doubling Property}

\subsection{Preparations for the Proof}
For the proof we need the following proposition of Xie and a result of Väisälä which we cite verbatim
\begin{proposition}[Proposition 3.6 in~\cite{xie_nagata_2008}]\label{prop:xie}
Let $f : (X_1, d_1) \to (X_2, d_2)$ be a quasi-Möbius homeomorphism. Then $f$ can be written as $f = f_2^{-1} \circ f' \circ f_1$, where
$f'$ is a quasi-symmetric map, and $f_i$ for $i \in \{1,2\}$ is either a metric inversion or the identity map on the metric space $(X_i,d_i)$.
\end{proposition}

\begin{proposition}[Theorem 3.10 in~\cite{vaisala_quasimobius_1984}]
Let $(X,d)$ be an unbounded metric space and let $f : X \to Y$ be a quasi-Möbius map. Then $f$ is quasi-symmetric
if and only if $f(x) \to \infty$ as $x \to \infty$.
If $X$ is any metric space and if $f: X \cup \{\infty\} \to Y \cup \{\infty\}$ is quasi-Möbius with $f(\infty) = \infty$, then
$f_{\vert X}$ is quasi-symmetric.
\end{proposition}

\begin{remark}
Let $(X,d)$ be an unbounded space. Then we can build the completed space with respect to the infinitely remote point
$\bar{X} := X \cup \{\infty\}$ together with an extended metric $\bar{d}$. Let $\bar{d}(x,y) := d(x,y)$ and $\bar{d}(\infty, x) :=  \bar{d}(x, \infty) = \infty$ for all $x,y \in X$. Furthermore
let $\bar{d}(\infty, \infty) = 0$. Then clearly $(X,d)$ is doubling if and only if $(\bar{X}, \bar{d})$ is doubling.
\end{remark}

\begin{theorem}\label{thm:inversion}
Let $(X,d)$ be an metric doubling space with doubling constant D, where $d$ is an extended metric~\cite{buyalo_elements_2007}
and denote by $\infty \in X$ the infinitely remote point in $(X,d)$.
Furthermore let $p \in X$ with $p \neq \infty$ and let $i_p$ be given by
$i_p(x,y) := \frac{d(x,y)}{d(p,x)\,d(p,y)}$ for all $x,y \in X\setminus\{\infty\}$ and $i_p(\infty,x) := i_p(x,\infty) := \frac{1}{d(p,x)}$.
Define $d_p(x,y) := \inf\{\sum_{i=1}^{k}i_p(x_i, x_{i-1}) \, \vert\, x= x_0, \ldots, x_k = y \in X\setminus\{p\}\}$.
Then $(X,d_p)$ is doubling with constant at most $D^{10}+1$.
\end{theorem}
\begin{proof}
  If $(X,d)$ is bounded, consider the space $(\bar{X}, \bar{d})$, with $\bar{X} := X \cup \{\infty\}$
  and $\bar{d}(x,y) := d(x,y)$ for all $x,y \in X$ and $\bar{d}(x, \infty) := \infty$.
  $(\bar{X},\bar{d})$ is doubling. Furthermore if $(\bar{X},\bar{d}_p)$ is doubling, then so is
  $(X, d_p)$. We therefore only need to show the theorem for unbounded $X$.

  We have the following relation for all $x,y \in X \setminus \{p\}$~\cite{buckley_metric_2008}:
  \[\frac{1}{4} i_p(x,y) \leq d_p(x,y) \leq i_p(x,y) \leq \frac{1}{d(x,p)} + \frac{1}{d(y,p)}.\]

Let $x_0 \in X\setminus\{p\}$ and $r >0$. Let
$B':= B'_r(x_0):=\{x \in X \,\vert\, d_p(x_0,x) \leq r\}$ be the ball of radius
$r$ in the space $(X,d_p)$. We consider the following two cases
\begin{enumerate}
\item If $B' \cap B'_{\frac{1}{2}r}(\infty) \neq \emptyset$, then let
  $A' := B'_r(x_0) \setminus B'_{\frac{1}{2}r}(\infty)$. Take $y_0 \in A'$. 
For any two points $x, y \in A'$ we have by definition of the metric $d_p$ and the above relation that
\[i_p(x,y) = \frac{d(x,y)}{d(p,x)\,d(p,y)} \leq 4 d_p(x,y) \leq 8 r,\]
and $\frac{1}{d(y,p)} = i_p(\infty,y) \geq d_p(\infty, y) > \frac{1}{2}r$. From this it follows that
\[d(x,y) \leq 8 r d(p,x)d(p,y) \leq \frac{32}{r}.\]
In particular we know that $A' \subseteq B_{\frac{32}{r}}(y_0) := \{x \in X \,\vert\, d(y_0,x) \leq \frac{32}{r}\}$. By the assumption we furthermore have for all $x \in B'$ that
\[d_p(x, \infty) \leq 2r + \frac{1}{2}r = \frac{5}{2}r\]
and therefore also
\[\frac{1}{d(p,x)} \leq \frac{5}{2}r,\] 
from which it follows that
\[d(p,x) \geq \frac{2}{5r}.\]
The space $(X,d)$ is doubling and we can find $D^N$ balls $b_i$ of radius $\frac{32}{r} 2^{-N}$ with centerpoints $x_i$ covering $B_{\frac{32}{r}}(y_0)$. Let
$\tilde{b}_i := b_i \cap A'$ then we have for all $x,y \in \tilde{b_i}$:
\[d_p(x,y) \leq i_p(x,y) = \frac{d(x,y)}{d(p,x)\,d(p,y)} \leq \frac{\frac{64}{r}\, 2^{-N}}{\frac{2}{5r}\frac{2}{5r}} = \frac{64 \cdot 5^2 \cdot r^2}{2^2\, 2^N\, r} 
= \frac{400}{2^N}\,r.\]
In particular for $N:= 10$ we know that we have constructed a cover of $B' \subseteq A' \cup B'_{\frac{1}{2}r}(\infty)$ by $D^{10}+1$ balls
of radius $\frac{1}{2}r$.
\item In case that $B' \cap B'_{\frac{1}{2}r}(\infty) = \emptyset$, we know that $d_p(x_0, \infty) > r$ and also
$d_p(B',\infty):= \inf_{x \in B'} d_p(x,\infty) \geq \frac{1}{2}r$. For all $y \in B'$ we have
\[i_p(x_0,y) = \frac{d(x_0,y)}{d(p,x_0)\,d(p,y)} \leq 4d_p(x_0,y) \leq 4r,\]
from which it follows that
\[d(x_0,y) \leq 4r d(p,x_0) d(p,y) \leq \frac{4r}{d_p(\infty, x_0)\, d_p(\infty,y)} \leq \frac{4r}{d_p{(\infty, B')}^2}.\]
We therefore have $B' \subseteq B_{\frac{4r}{d_p{(\infty, B')}^2}}(x_0)$ and by the doubling property of $(X,d)$ we can cover by $D^N$ balls $b_i$ of radius
$\frac{4r}{d_p{(\infty, B')}^2}\,2^{-N}$ with center points $x_i$. Let $\tilde{b}_i := b_i \cap B'$, then we have for any two
$x,y \in \tilde{b}_i$:
\[d_p(x,y) \leq i_p(x,y) = \frac{d(x,y)}{d(p,x)d(p,y)} \leq \frac{\frac{8r}{d_p{(\infty, B')}^2}\,2^{-N}}{d(p,x)\,d(p,y)}
= 2^{-N+4}\frac{d_p(\infty,x)\,d_p(\infty,y)}{d_p{(\infty, B')}^2}r.\]
Furthermore we have
\[d_p(x, \infty) \leq d_p(x_0,x) + d_p(x_0,\infty) \leq r + d_p(B', \infty) +r \leq 5 d_p(B', \infty).\]
In conclusion we get that
\[2^{-N+4}\frac{d_p(\infty,x)\,d_p(\infty,y)}{d_p{(\infty, B')}^2} \leq 2^{-N+4} \frac{5^2 d_p{(\infty, B')}^2}{d_p{(\infty, B')}^2} = \frac{8\cdot5^2}{2^N}.\]
It therefore follows that if we take $N:=9$, then we have a covering of $B'$ by $D^9$ balls of radius $\frac{1}{2}r$.
\end{enumerate}
\end{proof}

\begin{remark}
Note that if in addition $d \in \mathcal{M}$ where $(X, \mathcal{M})$ is Ptolemy Möbius, then $i_p = d_p$ and in particular
$(X,d_p)$ is doubling with constant at most $D^8 +1$.
\end{remark}

\subsection{Proof of \autoref{thm:main}}

\begin{proof}[Proof of Theorem~\ref{thm:main}]
It remains to show the theorem for $(X,d)$ being a doubling metric space, $f:(X,d) \to (X,d')$ a metric inversion
and we have the following cases to check:
\begin{enumerate}
 \item $(X,d)$ unbounded, $(X,d')$ bounded;\label{itm:prf1}
 \item $(X,d)$ and $(X,d)$ both unbounded but with different points at infinity.\label{itm:prf3}
\end{enumerate}
Case~\ref{itm:prf3} follows directly
from Theorem~\ref{thm:inversion}. In the situation of~\ref{itm:prf1}, $d'$ is a metric inversion $d_p$ where $p$ is an isolated point
in $X$. That is there exists a $\epsilon > 0$ such that $d(p,x) > \epsilon$ for all $x \in X \setminus \{p\}$. The proof of Theorem~\ref{thm:inversion} still holds.
\end{proof}

\section{Invariance of Uniform Disconnectedness}\label{sec:unidisc}

The proof of \autoref{thm:undisc} will again make use of some of the propositions from the previous sections.

\begin{figure}
  \begin{center}


\begin{tikzpicture}
      \foreach \a in {0,1,...,11}{
        \ifthenelse{\a=4 (\OR \a>8 \AND \a<9)}{
          \draw (\a*360/10+7*360/10: 3cm) node(x\a){};
        }{
          \draw (\a*360/10+7*360/10: 3cm) node(x\a){$\circ$};
        }
      }

      \draw [loosely dotted, thick] (3.5*360/10+7.33*360/10:2.9cm) -- (3.5*360/10+7.66*360/10:2.9cm);
      \draw [loosely dotted, thick] (3.5*360/10+11.75*360/10:2.9cm) -- (3.5*360/10+12.25*360/10:2.9cm);
  
      \node (p) at (0,0) {$\circ$} node [right] {$p$};

      \draw [dotted](x0) -- (x1) node [pos=0.5, below] {$l$};
      \draw [dotted](x1) -- (x2) node [pos=0.5, below right] {$l_0$};
      \draw [dotted](x2) -- (x3) node [pos=0.5, right] {$l_1$};
      \draw [dotted](p) -- (x0) node [pos=0.6] {$r_n$};
      \draw [dotted](p) -- (x1) node [pos=0.6] {$r_0$};
      \draw [dotted](p) -- (x2) node [pos=0.6] {$r_1$};

      \draw [dotted](p) -- (x6) node [pos=0.6] {$r_{s}$};
      \draw [dotted](p) -- (x7) node [pos=0.6] {$r_{s+1}$};
      \draw [dotted](x5) -- (x6) node [pos=0.5, above] {$l_{s-1}$};
      \draw [dotted](x6) -- (x7) node [pos=0.5, above left] {$l_s$};

      \draw (x0) node [below] {$x_n$};
      \draw (x1) node [below] {$x_0$};
      \draw (x2) node [below right] {$x_1$};
      \draw (x3) node [right] {$x_2$};
      \draw (x5) node [above] {$x_{s-1}$};
      \draw (x6) node [above] {$x_{s}$};
      \draw (x7) node [above left] {$x_{s+1}$};

      \draw [dotted] (x7)--(x8) node[pos=0.5, left] {$l_{s+1}$};
      \draw (x8) node [left] {$x_{s+2}$};
      
      \draw [dotted] (x9)--(x0) node[pos=0.5, below left] {$l_{n}$};
      \draw (x9) node [left] {$x_{n-1}$};

    \end{tikzpicture}
  \end{center}
  \caption{The view of the $\theta$-chain in $(X,d)$}\label{fig:chain}
\end{figure}
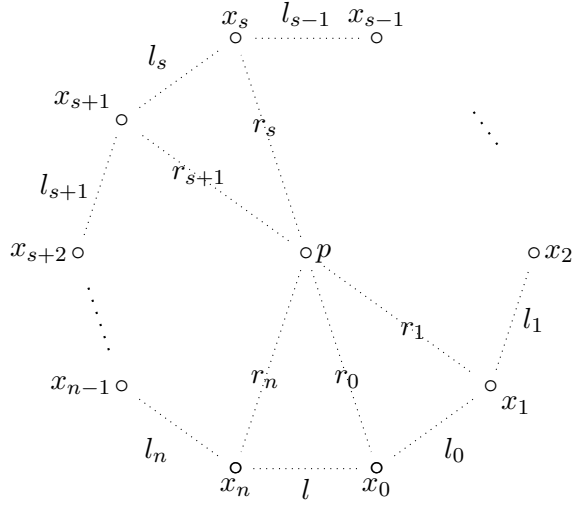

In the following let $(X,d)$ be a metric space, $p \in X$ and $\theta \leq \frac{1}{32}$.
We assume that $(X,d_p)$ is not $\theta$-uniformly disconnected, in particular there is some $\theta$-chain
$(x_0, x_1, \ldots, x_n)$ in $(X \setminus \{p\}, d_p)$. We keep this notation for the rest of this section.
In addition we introduce the following notation for convenience: Let $r_i := d(p,x_i)$, $l := d(x_0, x_n)$ and
$l_i := d(x_i, x_{i+1})$. This is illustrated in \autoref{fig:chain}.
Without loss of generality we can assume $r_n \geq r_0$.

\begin{remark}
  The condition for $(x_0, x_1, \ldots, x_n)$ being a $\theta$-chain in $(X,d_p)$ implies that
  \[\frac{l_i}{r_i r_{i+1}} \leq \frac{4 \theta l}{r_n r_0}\, \forall i \in \{0,\ldots, n-1\}.\]
  On the other hand if
  \[\frac{l_i}{r_i r_{i+1}} \leq \frac{\theta l}{4 r_n r_0}\, \forall i \in \{0,\ldots, n-1\}\]
  holds, then $(x_0, x_1, \ldots, x_n)$ is a $\theta$-chain in $(X,d_p)$.
\end{remark}

\begin{lemma}\label{lem:indexs}
  Assume that $(X,d)$ contains no $\sqrt[3]{4 \theta}$-chains. Then there is an index
  $s \in \{0, \ldots, n-1\}$ such that
  \[l_s > l \sqrt[3]{4 \theta}\]
  and
  \[\max\{r_s,r_{s+1}\} \sqrt[3]{4 \theta} \geq r_0.\]
\end{lemma}
\begin{proof}
  Assume for a contradiction that $r_s \sqrt[3]{4 \theta} < r_0$ and $r_{s+1} \sqrt[3]{4 \theta} < r_0$. Then
  from the condition in the remark above it follows
  \begin{equation}
    \frac{l_s}{r_s r_{s+1}} \leq \frac{4 \theta l}{r_n r_0} < \frac{4 \theta l_s}{\sqrt[3]{4 \theta} r_n r_0} <
    \frac{4 \theta l_s}{\sqrt[3]{4 \theta}^3 r_s r_{s+1}} =\frac{l_s}{r_s r_{s+1}}
  \end{equation}
  which is a contradiction.
\end{proof}

\begin{proposition}\label{prop:thetachain}
$(X,d)$ contains a $\sqrt[3]{4 \theta}$-chain.
\end{proposition}
\begin{proof}
  By the previous lemma we know that there must be some index $q$ such that $r_q \sqrt[3]{4 \theta} \geq r_0$
  and for all $i \in \{0, \ldots, q-1\}$ we have that $r_i \sqrt[3]{4 \theta} < r_0$.
  
  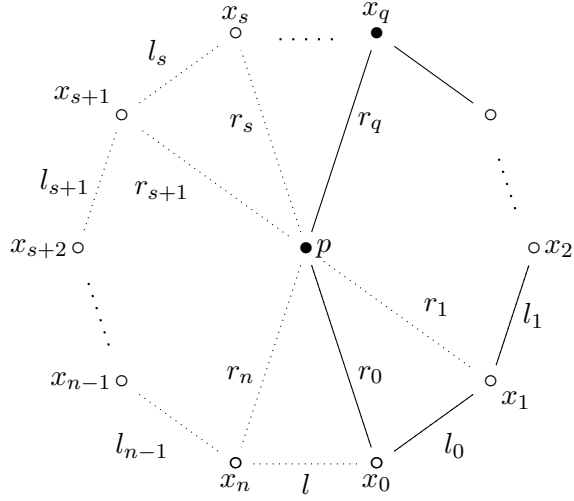
\begin{figure}
    \begin{center}


\begin{tikzpicture}
      \foreach \a in {0,1,...,11}{
        \ifthenelse{\a>8 \AND \a<9}{
          \draw (\a*360/10+7*360/10: 3cm) node(x\a){};
        }{
          \ifthenelse{\a = 5}{
            \draw (\a*360/10+7*360/10: 3cm) node(x\a){$\bullet$};
          }{
            \draw (\a*360/10+7*360/10: 3cm) node(x\a){$\circ$};
          }
        }
      }

      \draw [loosely dotted, thick] (3.5*360/10+8.8*360/10:2.8cm) -- (3.5*360/10+9.2*360/10:2.8cm);
      \draw [loosely dotted, thick] (3.5*360/10+6.8*360/10:2.8cm) -- (3.5*360/10+7.2*360/10:2.8cm);
      \draw [loosely dotted, thick] (3.5*360/10+11.75*360/10:2.9cm) -- (3.5*360/10+12.25*360/10:2.9cm);
  
      \node (p) at (0,0) {$\bullet$} node [right] {$p$};

      \draw [dotted](x0) -- (x1) node [pos=0.5, below] {$l$};
      \draw [dotted](x1) -- (x2) node [pos=0.5, below right] {$l_0$};
      \draw [dotted](x2) -- (x3) node [pos=0.5, right] {$l_1$};
      \draw [dotted](p) -- (x0) node [pos=0.6, left] {$r_n$};
      \draw [dotted](p) -- (x1) node [pos=0.6, right] {$r_0$};
      \draw [dotted](p) -- (x2) node [pos=0.6, above right] {$r_1$};

      \draw [dotted](p) -- (x6) node [pos=0.6, left] {$r_{s}$};
      \draw [dotted](p) -- (x7) node [pos=0.6, below left] {$r_{s+1}$};
      \draw [dotted](x6) -- (x7) node [pos=0.5, above left] {$l_s$};

      \draw (x0) node [below] {$x_n$};
      \draw (x1) node [below] {$x_0$};
      \draw (x2) node [below right] {$x_1$};
      \draw (x3) node [right] {$x_2$};
      \draw (x5) node [above] {$x_{q}$};
      \draw (x6) node [above] {$x_{s}$};
      \draw (x7) node [above left] {$x_{s+1}$};

      \draw [dotted] (x7)--(x8) node[pos=0.5, left] {$l_{s+1}$};
      \draw (x8) node [left] {$x_{s+2}$};
      
      \draw [dotted] (x9)--(x0) node[pos=0.5, below left] {$l_{n-1}$};
      \draw (x9) node [left] {$x_{n-1}$};

      \draw [] (p) -- (x5) node [pos=0.6, right] {$r_q$};
      \draw [] (x5) -- (x4);
      \draw [] (x3) -- (x2);
      \draw [] (x2) -- (x1);
      \draw [] (x1) -- (p);

    \end{tikzpicture}
    \end{center}
    \caption{The constructed $\sqrt[3]{4 \theta}$-chain in $(X,d)$}\label{fig:chain2}
  \end{figure}
  We claim that $(x_q, x_{q-1}, \ldots, x_1, x_0, p)$ is a $\sqrt[3]{4 \theta}$-chain in $(X,d)$. If this were not so,
  there would be some $i \in \{0, \ldots, q-1\}$ for which $r_q \sqrt[3]{4 \theta} < l_i$. But then
  \begin{equation}
    \frac{r_q \sqrt[3]{4 \theta}^2}{r_0 r_{q}} < \frac{r_q \sqrt[3]{4 \theta}}{r_i r_{q}}
    \leq \frac{r_q \sqrt[3]{4 \theta}}{r_i r_{i+1}} < \frac{l_i}{r_i r_{i+1}} \leq \frac{4 \theta l}{r_n r_0}
  \end{equation}
  implies
  \begin{equation}
    r_n < \sqrt[3]{4 \theta} l \leq \frac{1}{2} l
  \end{equation}
  which is a contradiction to the triangle inequality of the metric space $(X,d)$.
\end{proof}

\begin{proof}[Proof of Theorem~\ref{thm:undisc}]
  The proof of the theorem now follows directly from \autoref{prop:xie}.
\end{proof}

\section{Applications of the Theorems}

For the following we need a short definition~\cite{david_fractured_1997}:
Let $F$ be a finite set with $k \geq 2$ elements and let $F^\infty$ denote the set of sequences
$\{x_i\}_{i=1}^\infty$ with $x_i \in F$. For two elements $x = \{x_i\}, y = \{y_i\} \in F^\infty$
let \[L(x,y) = \sup\{I \in \mathbb{N}\,|\, \forall 1 \leq i \leq I: x_i = y_i\}.\]
In particular we have $L(x,x) = \infty$ and $L(x,y) = 0$ if $x_1 \neq y_1$.
Given $0 < a < 1$ set $\rho_a(x,y) = a^{L(x,y)}$. This defines an ultrametric on $F^\infty$.
We call $(F^\infty, \rho_a)$ the \emph{symbolic $k$-Cantor set with parameter $a$}.

As an application of the theorems we provide a generalization of the following result by David and Semmes:
\begin{proposition}[Proposition 15.11 (Uniformization) in \cite{david_fractured_1997}]
  Suppose that $(M, d)$ is a compact metric space which is bounded, complete, doubling, uniformly disconnected, and
  uniformly perfect. Then $M$ is quasi-symmetrically equivalent to the symbolic Cantor set $F^\infty$, where we
  take $F = \{0,1\}$ and we use the metric $\rho_a$ on $F^\infty$ with parameter $a = \frac{1}{2}$.
\end{proposition}

We can generalize this result as follows:

\begin{theorem}
  Suppose that $(M,d)$ is a complete, doubling, uniformly perfect and
  uniformly disconnected metric space. Then $M$ is quasi-Möbius
  equivalent to the symbolic Cantor set as given above.
\end{theorem}
\begin{proof}
  Let $p \in M$ be some point and let $s_p(x,y) = \frac{d(x,y)}{(d(x,p) + 1)(d(y,p) + 1)}$. 
  Let $\hat{d}_p(x,y) = \inf\{\sum_{i=1}^k s_p(x_i, x_{i-1}) : x = x_0,\ldots,x_k=y \in X\}$.
  We have~\cite{buckley_metric_2008}
  \[\frac{1}{4}s_p(x,y) \leq \hat{d}_p(x,y) \leq s_p(x,y) \leq \frac{1}{1+d(x,p)} + \frac{1}{1+d(y,p)}.\]
  Then
  the space $(M, \hat{d}_p)$ is bounded and satisfies all the properties of the above proposition:
  The map $f:(X,d) \to (X,\hat{d}_p)$ given by $d \mapsto \hat{d}_p$ is Möbius. By \autoref{thm:undisc} and
  \autoref{thm:main}, doubling and uniformly disconnectedness are invariant under Möbius maps.
  The invariance of uniformly perfectness follows from \cite{meyer_uniformly_2009}, and the
  invariance of completeness follows from \cite{beyrer_trees_2015}.
  Totally boundedness follows from the doubling property and therefore the space
  $(X,\hat{d}_p)$ is compact.
\end{proof}

We can apply the same idea to Proposition 16.9 in \cite{david_fractured_1997} and we get:
\begin{corollary}
  Let (M,d) be a complete Ahlfors regular metric space of dimension $\gamma$ which is uniformly disconnected.
  Then there exists a doubling measure $\mu$ on $F^\infty$, and $(M,d)$ is quasi-Möbius equivalent
  to $(F^\infty, D)$, where $D$ is given by
  \[D(x,y) = \left(\mu(\bar{B}(x,d_a(x,y))) + \mu(\bar{B}(y,d_a(x,y)))\right)^{\frac{1}{\gamma}},\]
  and $0 < a < 1$.
\end{corollary}

This follows from the above remarks and the invariance of Ahlfors regularity under $d \mapsto \bar{d}_p$
as shown in~\cite{li_preservation_2015}.

\newpage
\appendix
\section{Appendix}\label{apx:A}


\begin{proposition}
Let $(X,d)$ be a $K$-quasi-metric space~\cite{buyalo_elements_2007}. Let $X_\infty$ denote the infinite remote set and let
$\infty \in X_\infty$, i.e.\ the space satisfies the relations
\begin{enumerate}
 \item $d(x,y) = 0 \iff x = y$,
 \item $d(x,y) = d(y,x)$,
 \item $d(x,y) \leq K \max\{d(x,z), d(z,y)\}$ for all $x,y,z \in X$ for which all distances are defined,
 \item $d(x,y) < \infty \iff x,y \in X \setminus X_\infty$.
\end{enumerate}
Let $\lambda : X \to [0, \infty]$, $L > 0$ and $K' \geq K$ be such that $X_\infty = \lambda^{-1}(\infty)$ and
\begin{enumerate}
 \item $d(x,y) \leq K' \max\{L \lambda(x), L \lambda(y)\}$,
 \item $L \lambda(x) \leq K' \max \{d(x,y), L \lambda(y)\}$.
\end{enumerate}
Denote by $X'_\infty := \{\lambda^{-1}(0)\}$. Define a new metric $d_\lambda : (X \times X) \setminus (X'_\infty \times X'_\infty) \to [0, \infty]$ by
\begin{enumerate}
 \item $d_\lambda(x,y) := \frac{d(x,y)}{\lambda(x)\lambda(y)}$ for $x,y \in X \setminus X'_\infty$,
 \item $d_\lambda(x, \infty) := d_\lambda(\infty, x) := \frac{L}{\lambda(x)}$ for $\infty \in X_\infty$,
 \item $d_\lambda(\infty, \infty) = 0$ for $\infty \in X_\infty$,
 \item $d_\lambda(x, p) := d_\lambda(p,x) := \infty$ for $p \in X'_\infty$.
\end{enumerate}
If $(X, d)$ is doubling with constant $D$ then $(X, d_\lambda)$ is doubling with constant at most $D^{\lceil \log_2(8 K'^{10}K) \rceil} + 1$.
\end{proposition}
\begin{proof}
By Prop 5.3.6 in~\cite{buyalo_elements_2007}, $d_\lambda$ is a $K'^2$-quasi-metric. In particular we have for all $x,y,z \in X$ for which all distances
are defined, that:
\[d_\lambda(x,y) \leq K'^2 \max\{d(d,z), d(z,y)\}.\]
Let $x_0 \in X$, $x_0 \neq p \in X'_\infty$ and $r > 0$ and let $B' := B'_r(x_0) := \{x \in X \,\vert\, d_{\lambda}(x_0, x) \leq r\}$. Consider the following cases
\begin{enumerate}
 \item If $B' \cap B'_{\frac{1}{2}r}(\infty) \neq \emptyset$, then let $A' := B ' \setminus B'_{\frac{1}{2}r}(\infty)$. For all $x, y \in B'$ we have 
  \[d_\lambda(x,y) = \frac{d(x,y)}{\lambda(x)\lambda(y)} \leq K'^2 r,\]
  from which it follows that
  \[d(x,y) \leq K'^2r \lambda(x) \lambda(y).\]
 Furthermore we have for all $x \in A'$ that $d_\lambda(\infty, x) = \frac{L}{\lambda(x)} > \frac{1}{2}r$ and therefore also
 $\lambda(x) < \frac{2L}{r}$. Combining both equations we get that for all $x, y \in A'$ we have
  \[d(x,y) \leq K'^2r \frac{2L}{r} \frac{2L}{r} = \frac{K'^2 4L^2}{r}.\]
Without loss of generality assume $x_0 \in A'$. By the doubling property of $(X,d)$ we can cover $B_{\frac{K'^2 4L^2}{r}}(x_0)$ by at most $D^N$ balls $b_i$
of radius $\frac{K'^2 4L^2}{r} 2^{-N}$. Let $\tilde{b}_i := b_i \cap A'$ then we have for all $x,y \in \tilde{b}_i$:
\[d_\lambda(x,y) \leq \frac{\frac{K'^2 4L^2}{r2^{N}}K}{\lambda(x)\lambda(y)}.\]
 By the assumption there is a $\bar{x} \in B' \cap B'_{\frac{1}{2}r}(\infty)$ and we have for $x \in B'$ that $d_\lambda(x, \bar{x}) \leq K'^2 r$,
 therefore we also have $\frac{L}{\lambda(x)} = d_\lambda(x, \infty) \leq K'^4 r$ and $\lambda(x) \geq \frac{L}{K'^4r}$. In conclusion we get for all $x,y \in \tilde{b}_i$:
 \[d_\lambda(x,y) \leq \frac{\frac{K'^2 4L^2}{r2^{N}}K}{\lambda(x)\lambda(y)} \leq \frac{\frac{K'^2 4L^2}{r2^{N}}K}{\frac{L}{K'^4r} \frac{L}{K'^4r}}
=\frac{K'^{10}K4r}{2^N}.\]
In particular for $N := \lceil \log_2(8 K'^{10}K) \rceil$ we get a cover of $B'$ by at most $D^N +1$ balls of half the radius.
\item If $B' \cap B'_{\frac{1}{2}r}(\infty) = \emptyset$, then we have $d_\lambda(x_0, \infty) > r$ and $d_\lambda(B', \infty) > \frac{1}{2} r$. For all
$y \in B'$ we have $d_\lambda(x_0, y) = \frac{d(x_0, y)}{\lambda(x_0)\lambda(y)} \leq r$ and therefore also
\[d(x_0, y) \leq r \lambda(x_0) \lambda(y) \leq \frac{rL^2}{d_\lambda(\infty, x_0) d_\lambda(\infty, y)} = \frac{rL^2}{d_\lambda{(B', \infty)}^2}.\]
By the doubling property of $(X, d)$ we can find $D^N$ balls $b_i$ of radius $\frac{rL^2}{d_\lambda{(B', \infty)}^2} 2^{-N}$ covering
$B'$. Let $\tilde{b}_i := b_i \cap B'$, then we have for any $x,y \in \tilde{b}_i$:
\[d_\lambda(x,y) = \frac{d(x,y)}{\lambda(x)\lambda(y)} \leq \frac{K \frac{rL^2 2^{-N}}{d_\lambda{(B',\infty)}^2}}{\lambda(x)\lambda(y)}
= \frac{Kr2^{-N}d_\lambda(\infty, x)d_\lambda(\infty, y)}{d_\lambda{(B', \infty)}^2}.\]
Furthermore for any $x \in B'$ we have
 \[d_\lambda(x, \infty) \leq K'^2 \max\{d_\lambda(x_0,x),d_\lambda(x_0, \infty)\} \leq K'^2 r
\leq K'^2 2 d_\lambda(B', \infty).\]
We can combine the estimates to get
\[d_\lambda(x,y) \leq \frac{Kr2^{-N} K'^4 4 d_\lambda{(B', \infty)}^2}{d_\lambda{(B', \infty)}^2} = Kr2^{-N}K'^4 4.\]
In particular for $N:= \lceil \log_2(8 K K'^4) \rceil$ we have constructed a covering by $D^N$ balls of radius at most $\frac{1}{2}r$.
\end{enumerate}
\end{proof}

\begin{proposition}
Let $(X,d)$ be a $K$-quasi-metric space~\cite{buyalo_elements_2007}. Let $X_\infty$ denote the infinite remote set and let
$\infty \in X_\infty$, i.e.\ the space satisfies the relations
\begin{enumerate}
 \item $d(x,y) = 0 \iff x = y$,
 \item $d(x,y) = d(y,x)$,
 \item $d(x,y) \leq K \max\{d(x,z), d(z,y)\}$ for all $x,y,z \in X$ for which all distances are defined,
 \item $d(x,y) < \infty \iff x,y \in X \setminus X_\infty$.
\end{enumerate}
Let $\lambda : X \to [0, \infty]$, $L > 0$ and $K' \geq K$ be such that $X_\infty = \lambda^{-1}(\infty)$ and
\begin{enumerate}
 \item $d(x,y) \leq K' \max\{L \lambda(x), L \lambda(y)\}$,
 \item $L \lambda(x) \leq K' \max \{d(x,y), L \lambda(y)\}$.
\end{enumerate}
Denote by $X'_\infty := \{\lambda^{-1}(0)\}$. Define a new metric $d_\lambda : (X \times X) \setminus (X'_\infty \times X'_\infty) \to [0, \infty]$ by
\begin{enumerate}
 \item $d_\lambda(x,y) := \frac{d(x,y)}{\lambda(x)\lambda(y)}$ for $x,y \in X \setminus X'_\infty$,
 \item $d_\lambda(x, \infty) := d_\lambda(\infty, x) := \frac{L}{\lambda(x)}$ for $\infty \in X_\infty$,
 \item $d_\lambda(\infty, \infty) = 0$ for $\infty \in X_\infty$,
 \item $d_\lambda(x, p) := d_\lambda(p,x) := \infty$ for $p \in X'_\infty$.
\end{enumerate}
Let $\theta \leq \frac{1}{K^{19}}$. If $(X, d_\lambda)$ has a $\theta$-chain, then $(X,d)$ has a $\sqrt[3]{\theta K'^4}$-chain.
\end{proposition}
\begin{proof}
  Using the same notation as before in \autoref{sec:unidisc} we note that for all $i \in \{0, \ldots, n-1\}$ the
  following relation holds:
  \[\frac{l_i}{\frac{K'^2}{L^2}r_{i} r_{i+1}} \leq \frac{l_i}{\lambda(x_{i}) \lambda(x_{i+1})} \leq \frac{l \theta}{\lambda(x_0)\lambda(x_n)} \leq \frac{l \theta}{\frac{1}{K' L}r_0 r_n}.\]
  We can apply a similar argument as in \autoref{lem:indexs} to get an index $q$ for which
  \[r_0 \leq \sqrt[3]{\theta K'^4} r_q,\]
  and such that for all $i \in \{0, \ldots, q-1\}$ we have
  \[r_0 > \sqrt[3]{\theta K'^4} r_i.\]
  Assume again for a contradiction that $(x_q, x_{x-1}, \ldots, x_0, p)$ is not a $\sqrt[3]{\theta K'^4}$-chain. Then
  for some $i \in \{0, \ldots, q-1\}$:
  \begin{align}
    \frac{\sqrt[3]{\theta K'^4}^2r_q}{\frac{K'^2}{L^2} r_0 r_{q}}
    \leq\frac{\sqrt[3]{\theta K'^4}r_q}{\frac{K'^2}{L^2} r_i r_{q}}
    &\leq\frac{\sqrt[3]{\theta K'^4}r_q}{\frac{K'^2}{L^2} r_i r_{i+1}}
    \leq\frac{\sqrt[3]{\theta K'^4}r_q}{\lambda(x_i)\lambda(x_{i+1})}
    < \frac{l_i}{\lambda(x_i) \lambda(x_{i+1})}\\
    &\leq \frac{\theta l}{\lambda(x_0)\lambda(x_n)} 
    \leq \frac{\theta l}{\frac{1}{K'^2 L^2} r_0 r_n}
  \end{align}
  From this it follows that
  \[r_n < \sqrt[3]{\theta K'^4} K'^4 l \leq K^{-1} l.\]
\end{proof}


\bibliographystyle{plain}
\bibliography{../../../librarybt}

\bigskip\noindent
Loreno Heer\\
Institut für Mathematik\\
Universität Zürich\\
Winterthurerstrasse 190\\
CH-8057 Zürich\\
{\tt loreno.heer@math.uzh.ch}

\end{document}